\documentclass[oneside,12pt]{amsart}
\usepackage{hyperref}
\usepackage{amssymb,amscd,amsmath,latexsym,amsthm}
\usepackage[mathcal]{euscript}
\usepackage{xypic}
\xyoption{all}

\begin{document}

\title {Equivariant K-theory, groupoids and proper actions}

\author[Jose Cantarero]{Jose Cantarero}
\thanks{The author is partially supported by FEDER/MEC grant
MTM2010-20692.}
\thanks{Revised version published at Journal of K-theory: K-theory and its Applications to Algebra, Geometry and Topology, Volume 9, Issue 3 (2012), 475-501. \url{http://journals.cambridge.org/abstract_S1865243309998816}. Full copyright belongs to ISOPP.} 
\address{
Department of Mathematics\\
\newline \indent Stanford University\\
\newline \indent Stanford, California 94305\\
\newline \indent USA.}
\email{cantarer@stanford.edu}

\def \Pr{{\mathcal P}}
\def \P{{\mathbb P}}
\def \C{{\mathbb C}}            
\def \N{{\mathbb N}} 
\def \Z{{\mathbb Z}}
\def \KV{{\mathbb K}}
\def \R{{\mathbb R}}
\def \F{{\mathbb F}}
\def \G{{\mathcal G}}
\def \H{{\mathcal H}}
\def \K{{\mathcal K}}
\def \Fa{{\mathcal F}}
\def \Hy{{\mathbb H}}

\begin{abstract}
In this paper we define complex equivariant $K$-theory for actions of Lie groupoids using finite-dimensional vector
bundles. For a Bredon-compatible Lie groupoid $\G $, this defines a periodic cohomology theory on the category of finite 
$\G $-CW-complexes. We also establish an analogue of the completion theorem of Atiyah and Segal. Some examples are discussed.
\end{abstract}

\maketitle

\newtheorem{thm}{Theorem}[section]
\newtheorem{cor}[thm]{Corollary}
\newtheorem{lemma}[thm]{Lemma}
\newtheorem{prop}[thm]{Proposition}
\newtheorem{conj}[thm]{Conjecture}
\theoremstyle{definition}
\newtheorem{defn}[thm]{Definition}
\newtheorem{example}[thm]{Example}



\section{Introduction}

The recent theorem of Freed, Hopkins and Teleman \cite{FHT1, FHT2, FHT3, FHT4} relates the complex equivariant
twisted $K$-theory of a simply-connected compact Lie group acting on itself by conjugation to the Verlinde algebra. This result 
links information about the conjugation action with the action of the loop group on its universal space for proper actions. 
If we use the language of groupoids, the two associated groupoids are locally equivalent. The invariance of orbifold 
$K$-theory under Morita equivalence \cite{ALR} also seems to suggest that the language of groupoids is an appropriate 
framework to work with proper actions. 

The complex representation ring of a compact Lie group $G$ can be identified with the $G$-equivariant complex $K$-theory of a point. 
Equivariant complex $K$-theory is defined via equivariant complex bundles, but this procedure does not give a cohomology theory for proper 
actions of non-compact Lie groups in general, as shown in \cite{Ph}. Using infinite-dimensional complex $G$-Hilbert bundles, Phillips \cite{Ph}
constructs an equivariant cohomology theory for any second countable locally compact group $G$ on the category of proper locally compact $G$-spaces
that agrees with equivariant $K$-theory for actions of compact Lie groups. But in some cases it is enough to use finite-dimensional vector bundles. For example, L\"uck and Oliver show that they suffice for discrete groups \cite{LO}. Inspired by methods from that paper, we will construct a version of 
complex equivariant $K$-theory for actions of a Lie groupoid $\G = (\G_0,\G_1) $ (see the next section for some background on groupoids) by using 
extendable complex equivariant bundles, defined below. These bundles are finite-dimensional, but are required to satisfy an additional condition which 
will make sure we have a Mayer-Vietoris sequence. 

\begin{defn}

Let $X$ be a $\G$-space, $ \pi : X \longrightarrow \G_0 $ its anchor map and $ V \longrightarrow X $ a complex $\G$-vector bundle. We 
say $V$ is extendable if there is a complex $\G$-vector bundle $W \longrightarrow \G_0$ such that $ V $ is a direct summand of 
$ \pi ^* W $.

\end{defn}

The Grothendieck construction then gives a cohomology theory on the category of $\G$-spaces. For any $\G $-space $X$, 
$K_{\G}^*(X) $ is a module over $K_{\G}^*(\G_0)$ and the latter can be identified with $K_{orb}^*(\G)$ when $\G$ is an orbifold. 
But this theory does not necessarily satisfy Bott periodicity. In fact, it may not agree with classical equivariant $K$-theory 
when the action on the space is equivalent to the action of a compact Lie group. For our purposes we only need to study
$\G$-spaces that are built out of $\G$-cells, which are compact $\G$-spaces whose $\G$-action is equivalent to the action of a compact Lie 
group on a finite complex. 

\begin{defn}

A groupoid $\G$ is Bredon-compatible if given any $\G$-cell $U$, all $\G$-vector bundles on $U$ are extendable.

\end{defn}

Bredon-compatibility makes sure that $\G$-equivariant $K$-theory agrees with classical equivariant $K$-theory on the $\G$-cells. 
This condition also implies Bott periodicity for finite $\G$-CW-pairs, proved by a Mayer-Vietoris argument and induction over the 
cellular structure. These are two of the results in Section 3 needed to prove Theorem \ref{CohomologyTheory}, which
we reproduce here:   

\begin{thm}

If $\G $ is a Bredon-compatible Lie groupoid, the groups $K_{\G }^n(X,A) $ define a $ \Z /2 $-graded multiplicative 
cohomology theory on the category of finite $\G$-CW-pairs.  

\end{thm}

The purpose of this version of equivariant $K$-theory for actions of Lie groupoids is the existence of a completion theorem 
for proper actions of Lie groups. There are more general versions of equivariant $K$-theory for actions of locally compact groupoids 
available in the literature using $C^*$-algebras and $KK$-theory, such as the ones developed in \cite{EM} and \cite{LG}, but no completion 
theorems are expected to hold in such generality.

Emerson and Meyer consider the question of when equivariant $K$-theory for actions of locally compact groupoids on locally compact spaces
can be defined using finite-dimensional vector bundles in \cite{EM}. For instance, Theorem 6.4 in \cite{EM} shows that equivariant $K$-theory 
defined in terms of equivariant vector bundles is isomorphic to equivariant $K$-theory defined using $C^*$-algebras for actions of a second countable, 
locally compact, Hausdorff groupoid $\G$ with Haar system on a proper, $\G$-compact, second countable $\G$-space $X$ as long as the $C^*$-algebra of 
the action groupoid has an approximate unit of projections. This last condition holds if and only if for each $ x \in X $, and each irreducible 
representation $V$ of the stabilizer of $x$, there is a $\G$-vector bundle on $X$ whose fiber over $x$ contains $V$. We show in Section 3 that
equivariant $K$-theory defined using extendable vector bundles coincides with the versions defined in Definitions 2.2 and 2.3 of \cite{EM}. The 
following theorem corresponds to Theorem \ref{Cstar}:

\begin{thm}

Let $\G$ be a Bredon-compatible Lie groupoid such that for each $ x \in \G_0 $, and each irreducible representation $V$ of the stabilizer of $x$, 
there is a $\G$-vector bundle on $\G_0$ whose fiber over $x$ contains $V$. Then equivariant $K$-theory defined in terms of extendable $\G$-vector bundles 
coincides with equivariant $K$-theory defined in terms of $C^*$-algebras and $KK$-theory on the category of finite $\G$-CW-complexes.

\end{thm}

We introduce a universal $\G$-space $E\G$ as the limit of a sequence of free $\G$-spaces $E^n\G$ as in the case of compact Lie
groups. The quotient of $ E\G $ by the $\G $-action is $B\G $, the classifying space of $\G$. We form the fibered product 
$ X \times _{\pi } E^n\G $ over $\G_0$ and prove a generalization of the completion theorem of Atiyah and Segal \cite{AS} when 
$\G$ is finite.

\begin{defn}

A Lie groupoid $\G$ is finite if $\G_0$ is a finite $\G$-CW-complex and the spaces $ B^n\G = E^n\G / \G $ are compact.

\end{defn}

This version of $\G$-equivariant $K$-theory is invariant under weak equivalence. The naturality of completion maps and the Atiyah-Segal
completion theorem for actions of compact Lie groups imply the existence of a completion theorem for $\G$-cells. A spectral sequence argument
is then used to show that the completion map constructed at the end of Section 4 is an isomorphism. This result corresponds to Theorem 
\ref{Completion}, the main goal of Section 5:

\begin{thm}
Let $\G$ be a Bredon-compatible, finite Lie groupoid and $X$ a finite $\G$-CW-complex. Then we have an isomorphism of 
pro-rings 
\[ \{ K_{\G}^*(X)/I_{\G}^n K_{\G}^*(X) \} \longrightarrow \{ K_{\G}^*(X \times _{\pi} E^n\G ) \}. \]
In particular, for $X = \G_0$ we obtain an isomorphism of pro-rings
\[ \{ K_{\G}^*(\G_0)/I_{\G}^n K_{\G}^*(\G _0) \} \longrightarrow \{ K^*(B^n\G) \}. \]
\end{thm}

Finally, Section 6 studies the posibility of applying these result to proper actions of Lie groups which are not necessarily compact. These 
actions have an associated Lie groupoid which encodes the action of the group on its classifying space for proper actions. Actions of finite groups and 
compact Lie groups and proper actions of discrete groups, pro-discrete groups, almost compact groups and matrix groups give rise to Bredon-compatible 
groupoids. The finiteness condition is not automatic, but when it holds a completion theorem follows. 
\newline

\noindent {\bf Acknowledgements.} The author would like to thank the
referee, whose comments and suggestions really helped improve the 
quality and readibility of the paper.

\section{Background on groupoids}

In this section we review some basic facts about groupoids. All this material can be found in \cite{ALR} and \cite{M}.

\begin{defn} \label{Groupoid}
A topological groupoid $\G$ consists of a space $\G_0$ of objects and a space $\G_1$ of arrows, together with five continuous structure maps, listed below.

\begin{itemize}
\item The source map $s: \G_1 \rightarrow \G_0 $ assigns to each arrow $ g \in \G_1 $ its source $s(g)$.
\item The target map $t: \G_1 \rightarrow \G_0 $ assigns to each arrow $ g \in \G_1 $ its target $t(g)$. For two
objects $x$, $y \in \G_0$, one writes $ g: x \rightarrow y $ to indicate that $g \in \G_1$ is an arrow with $s(g)=x$
and $t(g)=y$.
\item If $g$ and $h$ are arrows with $s(h)=t(g)$, one can form their composition $hg$, with $s(hg)=s(g)$
and $t(hg)=t(h)$. The composition map $ m:\G_1 \times _{s,t} \G_1 \rightarrow \G_1 $, defined by $m(h,g)=hg$, is
thus defined on the fibered product
\[ \G_1 \times _{s,t} \G_1 = \{ (h,g) \in \G_1 \times \G_1 \mid s(h)=t(g) \} \]
and is required to be associative.
\item The unit map $ u:\G_0 \rightarrow \G_1 $ which is a two-sided unit for the composition. This means that $su(x) = x
= tu(x) $, and that $gu(x) = g = u(y)g $ for all $x$, $y \in \G_0$ and $ g: x \rightarrow y$.
\item An inverse map $ i: \G_1 \rightarrow \G_1 $, written $i(g) = g^{-1} $. Here, if $ g: x \rightarrow y$, then
$ g^{-1} : y \rightarrow x $ is a two-sided inverse for the composition, which means that $g^{-1}g = u(x) $ and 
$gg^{-1} = u(y)$.

\end{itemize}

\end{defn}

\begin{defn}

A Lie groupoid is a topological groupoid $\G$ for which $\G_0$ and $\G_1$ are smooth manifolds, and such that the structure
maps are smooth. Furthermore, $s$ and $t$ are required to be submersions so that the domain $\G_1 \times _{s,t} \G_1$ of $m$ 
is a smooth manifold.

\end{defn}

\begin{example}

Suppose a Lie group $K$ acts smoothly on a manifold $M$. One defines a Lie groupoid $K \rtimes M $ by $(K \rtimes M)_0 = M $
and $(K \rtimes M)_1 = K \times M $, with $s$ the projection and $t$ the action. Composition is defined from the multiplication
in the group $K$. This groupoid is called the action groupoid.

\end{example}

\begin{defn} \label{Orbits}

Let $\G $ be a Lie groupoid. For a point $x \in \G_0 $, the set of all arrows from $x$ to itself is a Lie group, denoted by
$\G _x$ and called the isotropy group at $x$. The set $ts^{-1}(x) $ of targets of arrows out of $x$ is called the orbit of $x$.
The quotient $ | \G | $ of $ \G_0 $ consisting of all the orbits in $ \G $ is called the orbit space. Conversely, we call
$ \G $ a groupoid presentation of $|\G | $.

\end{defn}

\begin{defn}

A Lie groupoid $\G$ is proper if $(s,t): \G_1 \rightarrow \G_0 \times \G_0 $ is a proper map. Note that in a proper Lie groupoid,
every isotropy group is compact.

\end{defn}

\begin{defn}

Let $\G $ and $\H $ be Lie groupoids. A strict homomorphism $\phi : \H \rightarrow \G $ consists of two smooth maps
$ \phi : \H_0 \rightarrow \G_0 $ and $ \phi : \H_1 \rightarrow \G_1 $ that commute with all the structure maps for the two
groupoids.

\end{defn}

\begin{defn}

A strict homomorphism $ \phi : \H \rightarrow \G $ between Lie groupoids is called an equivalence if:

\begin{itemize}
\item The map
\[ t\pi _1 : \G_1 \times _{s,\phi } \H_0 \rightarrow \G_0 \]
is a surjective submersion, where $ \G_1 \times _{s,\phi } \H_0 = \{ (g,y) \in \G_1 \times \H_0 \mid s(g) = \phi (y) \} $ and 
$ \pi _1 : \G_1 \times _{s,\phi } \H_0 \to \G_1$ is the projection to the first factor.
\item The square
\[ 
\diagram
\H_1 \rto^{\phi } \dto_{(s,t)} & \G_1 \dto_{(s,t)} \cr
\H_0 \times \H_0 \rto^{\phi \times \phi } & \G_0 \times \G_0
\enddiagram
\]
is a fibered product of manifolds.

\end{itemize}

\end{defn}

The first condition implies that every object $ x \in \G_0 $ can be connected by an arrow $ g: \phi (y) \rightarrow x $ to
an object in the image of $\phi $, that is, $\phi $ is essentially surjective as a functor. The second condition implies that
$\phi $ induces a diffeomorphism
\[ \H_1(y,z) \rightarrow \G_1(\phi (y), \phi (z)) \]
from the space of all arrows $ y \rightarrow z $ in $\H_1$ to the space of all arrows $ \phi (y) \rightarrow \phi (z) $ in
$\G_1$. In particular $\phi $ is full and faithful as a functor. 
\newline

A strict homomorphism $ \phi : \H \rightarrow \G $ induces a continuous map $ | \phi | : 
| \H | \rightarrow | \G | $. Moreover, if $\phi $ is an equivalence, $ | \phi | $ is a homeomorphism.

\begin{defn}

A local equivalence $ \phi : \H \longrightarrow \G $ is an equivalence with the additional property that each $ g_0 \in \G_0 $ has
a neighbourhood $U$ admitting a lift to $ \tilde{\H}_0 $ in the diagram

\[
\diagram
& \tilde{\H}_0 \dto \rto & \H_0 \dto_{\phi _0} \cr
& \G_1 \dto^s \rto^t & \G_0 \cr
U \morphism \dashed \tip \notip [-2,1] \rto & \G_0,
\enddiagram
\]

\noindent in which the square is a pullback square.

\end{defn}

\begin{defn} \label{C}

Two Lie groupoids $\G$ and $\G '$ are weakly equivalent if there exists a third groupoid $\H$ and two local equivalences
\[ \G \leftarrow \H \rightarrow \G '. \]

\end{defn}

\section{Groupoid actions and K-theory}

\begin{defn}

Let $\G$ be a groupoid. A (right) $\G$-space is a manifold $E$ equipped with an action of $\G$. Such an action is given
by two maps $\pi : E \rightarrow \G_0 $ (called the anchor map) and $ \mu : E \times _{\G_0} \G_1 \rightarrow E $. The latter map is 
defined on pairs $(e,g)$ with $ \pi (e) = t(g) $ and written $ \mu (e,g) = e \cdot g$. They must satisfy $\pi (e \cdot g)=s(g)$, 
$e \cdot 1_{\pi (e)} = e$ and $ (e \cdot g) \cdot h = e \cdot (gh) $. The space of orbits of this action is denoted by $ E/\G $.

\end{defn}

\begin{example}

If $\G$ is a groupoid, $ \G_1$ is a $\G$-space with the anchor map given by $ s: \G_1 \to \G_0 $, so that 
$ \G_1 \times_{\G_0} \G_1 = \G_1 \times_{s,t} \G_1 $. The action map is the composition $ m: \G_1 \times_{s,t} \G_1 \to \G_1 $
from Definition \ref{Groupoid}. Consider the following maps:
\newline
\newline
\noindent $ \alpha :\G_1 / \G \to \G _0 $
\newline  
\noindent $    \alpha([g])= t(g) $
\newline  
\newline
\noindent $ \beta : \G_0 \to \G_1 / \G $
\newline
\noindent $ \beta (x) = [u(x)] $,
\newline
\newline
where $[g]$ denotes the orbit of $ g \in \G _1 $. These maps define a homeomorphism between $ \G_1 /\G $ and $\G_0$. 
\end{example}

\begin{example}

Let $M$ be a $\G $-space. We can construct the action groupoid $ \H = \G \rtimes M $ which has space of objects $M$ and morphisms 
$ M \times _{\G_0} \G_1 $. This groupoid generalizes the earlier notion of action groupoid for a group action and the structure maps 
are formally the same as in that case.

\end{example}

\begin{defn}

Let $\G$ be a groupoid and let $X$, $Y$ be $\G$-spaces with anchor maps $\pi _X$ and $\pi _Y$, respectively. A map $ f : X \longrightarrow Y $ is 
$\G$-equivariant if it satisfies $ \pi _Y \circ f = \pi _X $ and $ f(x\cdot g) = f(x) \cdot g $ for all $ g \in \G _1 $ with $ t(g) = \pi _X (x) $.

\end{defn}

\begin{itemize}

\item $\G_0$ is a final object in the category of $\G $-spaces with the action given by $ e \cdot g = s(g) $ and anchor map given 
by the identity. The space of orbits for this action is $ | \G | $ (Definition \ref{Orbits}).

\item If $X$ and $Y$ are $\G $-spaces, the fibered product over $\G_0$,  $ X \times _{\pi } Y = \{ (x,y) \mid \pi (x) = \pi (y) \} $ becomes
a $\G $-space with coordinate-wise action. In particular $ X \times _{\pi } \G_0 = X $.

\item Similarly, if $X$ is a $\G$-space and $Y$ is any other space, $X \times Y$ is a $\G$-space with trivial action on the second
factor. In fact, $ X \times Y = X \times _{\pi } (Y \times \G_0) $.

\end{itemize}

\noindent All vector bundles in this paper are complex vector bundles.

\begin{defn}

Let $\G $ be a groupoid. A $\G $-vector bundle on a $\G $-space $X$ is a vector bundle $ p :E \longrightarrow X$ such that $E$ is a 
$\G $-space with fiberwise linear action and $p$ is a $\G$-equivariant map.

\end{defn}

\begin{defn}

Let $X$ be a $\G$-space and $ V \longrightarrow X $ a $\G$-vector bundle. We say that $V$ is $\G$-extendable if there is a $\G$-vector
bundle $W \longrightarrow \G_0$ such that $ V $ is a direct summand of $ \pi ^* W $.

\end{defn}

\begin{itemize}

\item Direct sum of $\G $-extendable vector bundles induces an operation on the set of isomorphism classes of  
$\G $-extendable vector bundles on $X$, making this set a monoid. We can also tensor $\G$-extendable vector bundles.

\item The pullback of a $\G$-extendable vector bundle by a $\G$-equivariant map is a $\G$-extendable vector bundle. 

\item All $\G$-vector bundles on $\G_0$ are $\G$-extendable. When $\G$ is an orbifold, $\G $-vector bundles on $\G_0$
correspond to $\G$-vector bundles over the orbifold $\G $ according to Definition 2.25 of \cite{ALR}.

\end{itemize}

\begin{defn}

Let $X$ be a $\G $-space. Then let Vect$_{\G } (X)$ be the set of isomorphism classes of extendable $\G $-vector bundles on $X$ and
$ K_{\G } (X) = K (\text{Vect}_{\G } (X)) $, where $K(A)$ denotes the Grothendieck group of a monoid $A$. We call $ K_{\G } (X) $ the 
extendable $\G $-equivariant $K$-theory of $X$.

\end{defn}

\noindent Let $X$ be a $\G$-space and $A$ a closed subspace of $X$ which is $\G$-invariant. We can now define the extendable $K$-groups as in 
Definition 3.1 of \cite{LO}:
\[ K^{-n}_{\G } (X) = \text{Ker} [ K_{\G }(X \times S^n) \stackrel{i^*}{\longrightarrow} K_{\G }(X) ], \]
\[ K^{-n}_{\G } (X,A) = \text{Ker} [ K^{-n}_{\G } ( X \cup _A X) \stackrel{j^*_2}{\longrightarrow} K^{-n}_{\G } (X) ], \] 

\noindent where $ i : X \rightarrow X \times S^n $ is the inclusion given by fixing a point in $ S^n $ and $ j_2 : X
\rightarrow X \cup _A X $ is one of the maps from $ X$ to the pushout. We equip $ X \times S^n $ with a $\G$-action by
taking as the anchor map the composition of the projection onto the first coordinate and the anchor map for $X$. Then we let
the groupoid act on the first coordinate. The anchor map for $ X \cup _A X $ is the only map to $\G_0$ making the pushout diagram 
commutative. The action is induced by the action of $\G$ on $X$.
\newline

From now on $\G$ will be a Lie groupoid. The following lemma follows easily from the definitions:

\begin{lemma} 

Let $(X,A)$ be a $\G$-pair. Suppose that $ X = \mathop{\amalg } _{i \in I} X_i $, the disjoint union of open $\G$-invariant subspaces 
$X_i$, and set $ A_i = A \cap X_i $. Then there is a natural isomorphism
\[  K_{\G}^{-n}(X,A) \longrightarrow \prod _{i \in I} K_{\G}^{-n}(X_i,A_i).  \]

\end{lemma}

\begin{defn}

A smooth left Haar system for a Lie groupoid $\G$ is a family $\{ \lambda ^a \mid a \in \G_0 \} $, where each $\lambda ^a$
is a positive, regular Borel measure on the manifold $ t^{-1}(a) $ such that:

\begin{itemize}

\item If $(V,\psi)$ is an open chart of $\G_1$ satisfying $ V \cong t(V) \times W $, and if $\lambda _W$ is the Lebesgue 
measure on $\R ^k$ restricted to $W$, then for each $ a \in t(V) $, the measure $ \lambda ^a \circ \psi $ is equivalent to 
$\lambda _W $, and the map $(a,w) \mapsto d(\lambda ^a \circ \psi _a )/d\lambda _W(w) $ belongs to 
$ C^{\infty} (t(V) \times W) $ and is strictly positive.

\item For any $ x \in \G_1 $ and $ f \in C^{\infty}_c(\G_1) $, we have
\[ \int _{t^{-1}(s(x))} f(xz)d\lambda ^{s(x)}(z) = \int _{t^{-1}(t(x))} f(y)d\lambda ^{t(x)}(y). \]

\end{itemize}

\end{defn}

\begin{prop}

Every Lie groupoid admits a smooth left Haar system.

\end{prop}

\noindent The proof can be found in \cite{P}, Theorem 2.3.1.
\newline

Note that we can use a smooth left Haar system to construct equivariant sections of $\G$-vector bundles from nonequivariant sections. More precisely, 
let $ f: X \to E $ be a nonequivariant section of a $\G$-vector bundle $ E \to X $. If we define 
$ h(x) = \int _{t^{-1}(\pi _X (x))} f(xg)g^{-1} d\lambda ^{\pi _X (x)}(g) $, then
\[ h(xk) = \int _{t^{-1}(\pi _X (xk))} f(xkg)g^{-1} d\lambda ^{\pi _X (xk)}(g) = \int _{t^{-1}(s(k))} f(xkg)g^{-1} d\lambda ^{s(k)}(g) = \]
\newline
\[ \int _{t^{-1}(s(k))} f(xkg)g^{-1}k^{-1} d\lambda ^{s(k)}(g) \cdot k = \int _{t^{-1}(t(z))} f(xz)z^{-1} d\lambda ^{t(z)}(z) \cdot k = \]
\newline
\[ \int _{t^{-1}(\pi _X (x))} f(xz)z^{-1} d\lambda ^{\pi _X(x)}(z) \cdot k = h(x) \cdot k. \]
\newline
\begin{cor}

If $f_0$, $f_1 : (X,A) \longrightarrow (Y,B) $ are $\G$-homotopic $\G$-maps between $\G$-pairs, then
\[  f_0^*=f_1^* : K_{\G}^{-n}(Y,B) \longrightarrow K_{\G}^{-n}(X,A)   \]
for all $n \geq 0 $.

\end{cor}

\begin{proof}

Since we can construct equivariant sections from nonequivariant sections, we can use the same argument from Propositions 1.1, 1.2 and 1.3 of \cite{SG3}. 
The relative case follows from the definition.
\end{proof}

The following lemma is a straightforward generalization of Lemma 1.5 in \cite{LO}:

\begin{lemma} \label{A}

Let $\phi : (X_1,X_0) \longrightarrow (X, X_2) $ be a map of $\G$-spaces, set $\phi _0 = \phi|_{X_0} $, and assume that 
$ X \cong X_2 \cup _{\phi _0} X_1 $. Let $p_1 : E_1 \longrightarrow X_1 $ and $p_2: E_2 \longrightarrow X_2$ be 
$\G$-extendable vector bundles, let $\overline{\phi }_0 : E_1|_{X_0} \longrightarrow E_2$ be an isomorphism of 
$\G$-vector bundles covering $ \phi _0 $, and set $ E = E_2 \cup _{\overline{\phi }_0} E_1 $. Then $ p = p_1 \cup p_2 : E \longrightarrow X $ is a 
$\G$-extendable vector bundle over $X$.

\end{lemma}

\begin{proof}

We have to show that $ p : E \to X $ is locally trivial. Since $ E_1 $ is locally trivial, so is $ E|_{X - X_2} \cong 
E|_{X_1 - X_0 } $. So it remains to find a neighbourhood of $X_2$ over which $E$ is locally trivial. Choose a closed neighbourhood
$W_1$ of $X_0$ in $X_1$ for which there is a strong deformation retraction $ r : W_1 \to X_0 $. By the homotopy invariance
for nonequivariant vector bundles over paracompact spaces, $r$ is covered by an isomorphism of vector bundles 
$\bar{r} : E_1|_{W_1} \to E_0 $ which extends $\bar{i}_1$. Set $W = X_2 \cup _{\phi _0} W_1 $. Then $\bar{r}$ extends, via
the pushout, to a map of vector bundles $E|_W \to E_2 $ which extends $ \bar{i}_2 $ and hence $ E|_W$ is locally trivial.
\end{proof}

The next lemma is a fundamental piece in the proof of the existence of a Mayer-Vietoris long exact sequence for $\G$-equivariant $K$-theory. It
was inspired by Lemma 3.7 of \cite{LO}. 

\begin{lemma}

Let $\phi : X \longrightarrow Y $ be a $\G$-equivariant map and let $E' \longrightarrow X $ be a $\G$-extendable vector bundle. 
Then there is a $\G$-extendable vector bundle $ E \longrightarrow Y $ such that $E'$ is a summand of $ \phi ^*E $.

\end{lemma}

\begin{proof} 

Consider $\pi : Y \longrightarrow \G_0 $ and $ \pi \phi : X \longrightarrow \G_0 $. Since
$E'$ is extendable, there is a $\G$-vector bundle $V$ on $\G_0$ such that $E'$ is a direct summand of $ (\pi \phi )^*V $.
Let $ E = \pi ^* V $. Then $E$ is a $\G$-vector bundle on $Y$ and it is the pullback of an extendable $\G$-vector bundle, hence it is
extendable. And we have that $ E' $ is a direct summand of $ (\pi \phi )^*V = \phi ^* E $.
\end{proof}

The existence of a Mayer-Vietoris sequence is shown next. The proof is an adaptation of Lemma 3.8 in \cite{LO}.

\begin{lemma}

Let
\[
\diagram
A \rto^{i_1} \dto_{i_2} & X_1 \dto_{j_1} \cr
X_2 \rto^{j_2} & X \\
\enddiagram
\]
be a pushout square of $\G$-spaces, where $i_1$ and $i_2$ are cofibrations. Then there is a natural exact sequence, infinite to the left
\[ \label{E1} \ldots \stackrel{d^{-n-1}}{\longrightarrow} K_{\G}^{-n}(X) \stackrel{j_1^* \oplus j_2^*}{\longrightarrow} K_{\G}^{-n}(X_1) \oplus
    K_{\G}^{-n}(X_2) \stackrel{i_1^*-i_2^*}{\longrightarrow} K_{\G}^{-n}(A) \stackrel{d^{-n}}{\longrightarrow} \ldots \]
\[  \ldots \longrightarrow K_{\G}^{-1}(A) \stackrel{d^{-1}}{\longrightarrow} K_{\G}^0(X) \stackrel{j_1^* \oplus j_2^*}{\longrightarrow} 
    K_{\G}^0(X_1) \oplus K_{\G}^0(X_2) \stackrel{i_1^*-i_2^*}{\longrightarrow} K_{\G}^0(A). \] 

\end{lemma}

\begin{proof}

We first show that the sequence
\begin{equation} \label{E2}
K_{\G}(X) \stackrel{j_1^* \oplus j_2^*}{\longrightarrow} K_{\G}(X_1) \oplus K_{\G}(X_2) \stackrel{i_1^* - i_2^*}{\longrightarrow} K_{\G}(A) 
\end{equation}
is exact, and hence the long sequence in the statement of the lemma is exact at $ K_{\G}^{-n}(X_1) \oplus K_{\G}^{-n}(X_2) $
for all $n$. Clearly the composite is zero. So fix an element $(\alpha _1, \alpha _2) \in \text{Ker}(i_1^*-i_2^*) $. By the previous
lemma, we can add an element of the form $ ([j_1^*E'], [j_1^*E']) $ for some $\G$-vector bundle $ E' \rightarrow X $, and
arrange that $ \alpha _1 = [E_1] $ and $ \alpha _2 = [E_2] $ for some pair of $\G$-vector bundles $ E_k \rightarrow X_k $. Then
$i_1^*E_1 $ and $ i_2^*E_2 $ are stably isomorphic, and after adding the restrictions of another bundle over $X$, we can
arrange that $ i_1^*E_1 \cong i_2^*E_2 $. Lemma \ref{A} now applies to show that there is a $\G$-vector bundle $E$ over $X$ such
that $j_k^*E \cong E_k $ for $ k = 1,2 $, and hence that $ (\alpha _1, \alpha _2) = ([E_1],[E_2]) \in \text{Im}(j_1^* \oplus j_2^*) $.
\newline

Assume now that $A$ is a retract of $X_1$. We claim that in this case the map
\begin{equation} \label{E3}
 \text{Ker} [K_{\G}(X) \stackrel{j_2^*}{\longrightarrow} K_{\G}(X_2)] \stackrel{j_1^*}{\longrightarrow} \text{Ker} [K_{\G}(X_1) \stackrel{i_1^*}{\longrightarrow} 
K_{\G}(A)] 
\end{equation}
is an isomorphism. It is surjective by the exactness of \eqref{E2}. So fix an element $[E]-[E'] \in \text{Ker}(j_1^* \oplus j_2^*) $. To
simplify the notation, we write $ E|_A = i_2^*j_2^*E $ and $ E|_{X_i} = j_i^*E$ for $i =1,2$. Let $p_1 : X_1 \to A $ be a retraction,
and let $ p : X \to X_2 $ be its extension to $X$. By the previous lemma, we can arrange that $ E|_{X_k} \cong E'|_{X_k} $
for $ k = 1,2 $. Applying the same lemma to the retraction $ p : X \to X_2 $, we obtain a $\G$-vector bundle $ F' \to X_2 $
such that $E'$ is a summand of $p^*F'$. Stabilizing again, we can assume that $E' \cong p^*F' $ and hence that $ F' \cong
E'|_{X_2} $ and $ E'|_{X_1} \cong p_1^*(F'|_A) \cong p_1^*(E'|_A) $. Fix isomorphisms $ \psi _k : E|_{X_k} \to E'|_{X_k} $
covering the identity on $X$. The automorphism $ (\psi |_A) \circ (\psi _1 |_A)^{-1} $ of $ E'|_A $ pulls back, under $p_1$,
to an automorphism $\phi $ of $ E'|_{X_1} $. By replacing $\psi _1 $ by $ \phi \circ \psi _1 $ we can arrange that 
$ \psi _1 |_A = \psi _2|_A $. Then $ \psi _1 \cup \psi _2 $ is an isomorphism from $E$ to $E'$, and this proves the exactness.
\newline

Now for each $ n \geq 1 $,
\[ K_{\G}^{-n}(A) = \text{Ker} [K_{\G}(A \times S^n) \to K_{\G}(A)] \]
\[ \cong \text{Ker} [K_{\G}(X \cup _{A \times pt} (A \times S^n)) \stackrel{\text{incl}^*}{\longrightarrow} K_{\G}(X)] \]
\[ \cong \text{Ker} [K_{\G}((X_1 \times D^n) \cup _{A \times S^{n-1}} (X_2 \times D^n)) \stackrel{(-,pt)^*}{\longrightarrow} K_{\G}(X)], \]
the last step since $(X_1 \times pt \cup A \times D^n )$  is a strong deformation retract of $ X_1 \times D^n $. Denote 
$ Y = (X_1 \times D^n) \cup _{A \times S^{n-1}} (X_2 \times D^n) $ and define
$ d^{-n} : K_{\G}^{-n}(A) \to K_{\G}^{-n+1}(X) $ to be the homomorphism which makes the following diagram commute:
\[
\diagram
0 \rto & K_{\G}^{-n}(A) \rto \dto^{d^{-n}} & K_{\G}(Y) \rto^{(-,pt)^*} \dto^{\text{incl}^*} & K_{\G}(X) \rto \dto^{Id} & 0 \\
0 \rto & K_{\G}^{-n+1}(X) \rto & K_{\G}(X \times S^{n-1}) \rto^{\hspace{0.75cm}(-,pt)^*} & K_{\G}(X) \rto & 0. \\
\enddiagram
\]

We have already shown that the desired long sequence is exact at $K_{\G}^{-n}(X_1) \oplus K_{\G}^{-n}(X_2) $ for all $n$. Let us
denote $Z = (X_1 \times D^n) \amalg (X_2 \times D^n) $ and \mbox{$ W = (X_1 \times D^n) \cup _{A \times pt} (X_2 \times D^n) $.} To prove exactness
at $ K_{\G}^{-n+1}(X)$ and $K_{\G}^{-n}(A) $ for any $ n \geq 1 $, apply the exactness of \eqref{E3} to the following split inclusion
of pushout squares:
\[
\diagram
X_1 \amalg X_2 \rto \dto & X_1 \amalg X_2 \dto  & (X_1 \amalg X_2) \times S^{n-1} \rto \dto & Z \dto \cr
X \rto \dto &       X  \rto^{\text{incl}} \dto & X \times S^{n-1} \rto \dto &  Y \dto \cr
X \rto &       X         & X \rto &  W. \\
\enddiagram
\]

The upper pair of squares induces a split surjection of exact sequence whose kernel yields the exactness of the long sequence at $K_{\G}^{-n+1}(X) $.
And since
\[ \text{Ker}[K_{\G}(W) \to K_{\G}(X)] \cong \text{Ker} \Big{[} K_{\G}(Z) \to K_{\G}(X_1 \amalg X_2) \Big{]} \cong K_{\G}^{-n}(X_1) \oplus K_{\G}^{-n}(X_2) \]
by \eqref{E3}, the lower pair of squares induces a split surjection of exact sequences whose kernel yields the exactness of the long sequence at 
$K_{\G}^{-n}(A)$.
\end{proof}

We now consider products on $K_{\G}^*(X) $ and on $ K_{\G}^*(X,A) $. We follow the analogous construction from Section 3 of \cite{LO}. 
Tensor products of $\G$-extendable vector bundles makes $K_{\G}(X) $ into a commutative ring, and all induced maps 
$ f^* : K_{\G}(Y) \longrightarrow K_{\G}(X) $ are ring homomorphisms. For each $ n,m \geq 0 $,
\[ K_{\G}^{-n-m}(X) \cong \text{Ker} [ K_{\G}^{-m}(X \times S^n) \longrightarrow K_{\G}^{-m}(X) ] = \]
\[ = \text{Ker} [ K_{\G}(X \times S^n \times S^m) \longrightarrow K_{\G}(X \times S^n) \oplus K_{\G}(X \times S^m)], \] 
where the first isomorphism follows from the usual Mayer-Vietoris sequences. Hence the composition
\[ K_{\G}(X \times S^n) \otimes K_{\G}(X \times S^m) \stackrel{p_1^* \otimes p_2^*}{\longrightarrow} 
K_{\G}(X \times X \times S^n \times S^m) \longrightarrow K_{\G}(X \times S^n \times S^m) \]
restricts to a homomorphism
\[ K_{\G}^{-n}(X) \otimes K_{\G}^{-m}(X) \longrightarrow K_{\G}^{-n-m}(X). \]

By applying the above definition with $ n = 0 $ or $m = 0 $, the multiplicative identity for $K_{\G}(X) $ is seen to be an 
identity for $K_{\G}^*(X)$. Associativity of the graded product is clear and graded commutativity follows upon showing that 
composition with a degree $-1$ map $ S^n \rightarrow S^m $ induces multiplication by $ -1 $ on $ K^{-n}(X)$. This product makes 
$K_{\G}^*(X)$ into a graded ring. Clearly, $f^* : K_{\G}^*(Y) \longrightarrow K_{\G}^*(X)$ is a ring homomorphism for any 
$\G$-map $ f: X \longrightarrow Y $. This makes $K_{\G}^*(X)$ into a $K_{\G}^*(\G_0)$-algebra, since $\G_0$ is a final object
in the category of $\G$-spaces.

It remains to prove Bott periodicity. Recall that $ \tilde{K}(S^2) = \text{Ker} [ K(S^2) \longrightarrow K(pt) ] \cong \Z $, and is 
generated by the Bott element $ B \in \tilde{K}(S^2) $, the element $ [ S^2 \times \C ] - [H] \in \tilde{K}(S^2) $, where $H$ 
is the canonical complex line bundle over $ S^2 = \C P^1 $. For any $\G$-space $X$, there is an obvious pairing
\[ K_{\G}^{-n}(X) \otimes \tilde{K}(S^2) \stackrel{\otimes}{\longrightarrow} \text{Ker} [ K_{\G}^{-n}(X \times S^2) \longrightarrow 
K_{\G}^{-n}(X) ] \cong K_{\G}^{-n-2}(X) \]
induced by external tensor product of bundles. Evaluation at the Bott element now defines a homomorphism
\[ b = b(X) : K_{\G}^{-n}(X) \longrightarrow K_{\G}^{-n-2}(X), \]
which by construction is natural in $X$. And this extends to a homomorphism
\[ b = b(X,A) : K_{\G}^{-n}(X,A) \longrightarrow K_{\G}^{-n-2}(X,A) \]
defined for any $\G$-pair $(X,A)$ and all $ n \geq 0 $.

\begin{defn}

An $n$-dimensional $\G$-cell is a space of the form $D^n \times U$ where $U$ is a compact $\G$-space such that $ \G \rtimes U $ is weakly
equivalent to an action groupoid corresponding to an action of a compact Lie group $G$ on a finite $G$-CW-complex. 

\end{defn}

\begin{defn} \label{B}

A $\G$-CW-complex $X$ is a $\G$-space together with an $\G$-invariant filtration
\[ \emptyset = X_{-1} \subseteq X_0 \subseteq X_1 \subseteq \ldots \subseteq X_n \subseteq \ldots \subseteq \mathop {\cup } 
_{n \geq 0} X_n = X \]
such that $X_n$ is obtained from $X_{n-1}$ for each $n \geq 0$ by attaching equivariant $n$-dimensional cells and $X$ carries the colimit topology with 
respect to this filtration. That is, for each $n \geq 0 $ there is a collection of $n$-cells $ \{ U_i \times D^n \mid i \in I_n \} $ for some index set 
$I_n$, along with attaching maps $q_i^n : U_i \times S^{n-1} \to X_{n-1}$, $Q_i^n : U_i \times D^n \to X_n $, and a $\G$-pushout
\[ \begin{CD}
\mathop{\coprod } \limits_{i \in I_n} U_i \times S^{n-1} @>\mathop{\coprod } \limits_{i \in I_n} q_i^n>> X_{n-1} \\
@VVV @VVV \\
\mathop{\coprod } \limits_{i \in I_n} U_i \times D^n @>\mathop{\coprod } \limits_{i \in I_n} Q_i^n>> X_n .
\end{CD} \]
Furthermore, we say that $X$ is a finite $\G$-CW-complex if it is constructed with a finite number of $\G$-cells, and similarly, 
$(X,A)$ is a finite $\G$-CW-pair if $X$ can be constructed from $A$ by attaching a finite number of $\G$-cells.
\end{defn}

\begin{defn}

A groupoid $\G$ is Bredon-compatible if given any $\G$-cell $U$, all $\G$-vector bundles on $U$ are extendable.

\end{defn}

\begin{example}

An example of a Bredon-compatible groupoid is $ \G = G \rtimes M $, where $G$ is a compact Lie group and $M$ is a finite
$G$-CW-complex. A $\G$-cell $U$ is a finite $G$-CW-complex with an equivariant map to $M$ and $\G$-vector bundles on $U$ are just $G$-vector
bundles. By Proposition 2.4 in \cite{SG3}, for any $G$-vector bundle $A$ on $U$, there is another $G$-vector bundle $B$ such that
$ A \oplus B $ is a trivial bundle, that is, the pullback of a $G$-vector bundle $V$ over a point. Consider the only map from $M$
to a point. The pullback of $V$ over this map is a $G$-vector bundle on $M$. If we pull it back to $U$ we recover
$ A \oplus B $ and therefore $\G$ is Bredon-compatible.

\end{example}

\begin{prop}

Suppose that $ F : \G \longrightarrow \H $ is a local equivalence. Then the pullback functor
\[ F^* : \{ \text{Fiber bundles on } \H \} \longrightarrow \{ \text{Fiber bundles on } \G \} \]
is an equivalence of categories.

\end{prop}

\begin{proof}

See Proposition A.18 of \cite{FHT2}. 
\end{proof}

\begin{cor}

If $\G$ is Bredon-compatible and $U$ is a $\G$-cell, then $K_{\G}^*(U) \cong K_G^*(M) $ for some compact Lie group
$G$ and a finite $G$-CW-complex $M$.

\end{cor}

\begin{thm} \label{Bott}

If $\G$ is Bredon-compatible, the Bott homomorphism
\[ b = b(X,A) : K_{\G}^{-n}(X,A) \longrightarrow K_{\G}^{-n-2}(X,A) \]
is an isomorphism for any finite $\G$-CW-pair $(X,A)$ and all $n \geq 0 $.

\end{thm}

\begin{proof}

Assume first that $ X = Y \cup _{\phi} (U \times D^m) $ where $U \times D^m$ is a $\G$-cell. Assume inductively that $b(Y)$ is
an isomorphism. Since $ K_{\G}^{-n}(U \times S^{m-1}) \cong K_G^{-n}( M \times S^{m-1}) $ and 
$ K_{\G}^{-n}(U \times D^m) \cong K_G^{-n}( M \times D^m) $, the Bott homomorphisms $b(U \times S^{m-1})$ and 
$b(U \times D^m)$ are isomorphisms by the equivariant Bott periodicity theorem for actions of compact Lie groups. The Bott map
is natural and compatible with the boundary operators in the Mayer-Vietoris sequence for $Y$, $X$, $U \times S^{m-1}$ and
$ U \times D^m $ and so $b(X)$ is an isomorphism by the $5$-lemma. The proof that $b(X,A)$ is an isomorphism follows immediately from the
definitions of the relative groups.
\end{proof}

\noindent Based on the Bott isomorphism we just proved, we can now redefine for all $n \in \Z $
\[ K_{\G}^n(X,A) = \left \{ \begin{array}{cc}
K_{\G}^0(X,A) & \text{if $n$ is even} \\
K_{\G}^{-1}(X,A) & \text{if $n$ is odd.} \end{array} \right. \]
For any finite $\G$-CW-pair $(X,A)$, define the boundary operator $\delta ^n : K_{\G}^n(A) \longrightarrow K_{\G}^{n+1}(X,A) $ to
be $ \delta : K_{\G}^{-1}(A) \longrightarrow K_{\G}^0(X,A) $ if $n$ is odd, and to be the composite
\[ K_{\G}^0(A) \stackrel{b}{\longrightarrow} K_{\G}^{-2}(A) \stackrel{\delta ^{-2}}{\longrightarrow} K_{\G}^{-1}(X,A) \]
if $n$ is even.

We collect now all the information we have about $\G$-equivariant $K$-theory in the following theorem:

\begin{thm} \label{CohomologyTheory}

If $\G $ is a Bredon-compatible Lie groupoid, the groups $K_{\G }^n(X,A) $ define a $ \Z /2 $-graded multiplicative cohomology 
theory on the category of finite $\G$-CW-pairs.  

\end{thm}

There are more general versions of equivariant $K$-theory for actions of locally compact groupoids 
available in the literature using $C^*$-algebras and $KK$-theory, such as the ones developed in 
\cite{EM} and \cite{LG}, but no completion theorems are expected to hold for them. For example, Emerson
and Meyer define two versions of equivariant $K$-theory for actions of a second countable, locally compact, Hausdorff 
groupoid $\G$ in \cite{EM}. The first definition uses the $C^*$-algebra of the action groupoid and it is given by
$CK_{\G}^*(X) = K_*(C^*(\G \rtimes X)) $. We use this notation to avoid confusion with our version of $\G$-equivariant
$K$-theory. The second definition is a Kasparov type $K$-theory, called $\G$-equivariant representable $K$-theory and given
by $RK_{\G}^*(X) = KK_*^{\G \rtimes X}(C_0(X),C_0(X))$. These definitions are extended to pairs and shown to be cohomology
theories on a suitable category of $\G$-spaces. They are isomorphic if $X$ is $\G$-compact, that is, if $X/\G$
is compact. 

Emerson and Meyer consider the question of when equivariant $K$-theory for groupoid actions can be defined using finite-dimensional
vector bundles. Let $\G$ be a second countable, locally compact, Hausdorff groupoid with Haar system, $X$ a proper, $\G$-compact, 
second countable $\G$-space and $A$ a closed $\G$-invariant subset of $X$. Assume that for each $ x \in X $, and each irreducible representation $V$ of 
the stabilizer of $x$, there is a $\G$-vector bundle on $X$ whose fiber over $x$ contains $V$. Note that since the anchor map $ \pi : X \to \G _0 $ is 
equivariant, it is sufficient that this condition holds for the action of $\G$ on $\G_0$. Then Theorems 6.4 and 6.14 in \cite{EM} show that equivariant 
$K$-theory defined in terms of $\G$-vector bundles $VK_{\G}^*(X,A)$ is isomorphic to $CK_{\G}^*(X,A)$ and $RK^*_{\G}(X,A)$.  

If moreover, $\G$ is a Bredon-compatible Lie groupoid and $U$ is a $\G$-cell, all $\G$-vector bundles on $U$ are extendable and therefore
$K_{\G}^*(U) \cong VK_{\G}^*(U) $. Assume that for each $ x \in \G_0 $, and each irreducible representation $V$ of the stabilizer of $x$, 
there is a $\G$-vector bundle on $\G_0$ whose fiber over $x$ contains $V$, so that $ VK_{\G}^*(Y) = CK_{\G}^*(Y)$ for any proper, $\G$-compact, second countable $\G$-space $Y$. The same argument used in the proof of Theorem \ref{Bott} implies that $K_{\G}^*(X) \cong CK_{\G}^*(X) 
\cong RK_{\G}^*(X) $ for all finite $\G$-CW-complexes. Therefore we have proved the following theorem:

\begin{thm} \label{Cstar}

Let $\G$ be a Bredon-compatible Lie groupoid such that for each $ x \in \G_0 $, and each irreducible representation $V$ of the stabilizer of $x$, 
there is a $\G$-vector bundle on $\G_0$ whose fiber over $x$ contains $V$. Then equivariant $K$-theory defined in terms of extendable $\G$-vector bundles 
coincides with equivariant $K$-theory defined in terms of $C^*$-algebras and $KK$-theory (\cite{EM}, Definitions 2.2 and 2.3) on the category of 
finite $\G$-CW-complexes.

\end{thm}

The condition on $\G_0$ is satisfied, for instance, by orbifolds (Example 6.17 in \cite{EM}) and by the action groupoids associated to certain actions of 
locally compact groups $G$ on proper $G$-compact spaces (Theorem 6.15 in \cite{EM}).

\section{The completion map}

Given a Lie groupoid $\G$, we can associate an important topological space to it, namely its classifying space 
$B \G $. For $ n \geq 1 $, let $\G_n$ be the iterated fibered product
\[ \G_n = \{ (g_1,\ldots ,g_n) \mid g_i \in \G_1, s(g_i)=t(g_{i+1}) \text{ for }i=1,\ldots ,n-1 \}. \]
Together with the space of objects $\G_0$, the spaces $\G_n$ have a simplicial structure called the nerve of $\G$. Here we are
really thinking of $\G$ as a topological category. Following the usual convention, we define face operators $d_i : \G_n \rightarrow
\G_{n-1} $ for $ i = 0,...,n$, given by
\[ d_i (g_1,\ldots ,g_n) = \left\{ \begin{array}{lll}
                               (g_2,\ldots ,g_n) & \mbox{if $i=0$}\\
                               (g_1,\ldots ,g_{n-1}) & \mbox{if $i=n$}\\
                               (g_1,\ldots ,g_ig_{i+1},\ldots ,g_n) & \mbox{otherwise}\end{array} \right. \]
for $ 0 < i < n $ when $ n > 1 $. Similarly, we define $ d_0(g) = s(g) $ and $ d_1(g) = t(g) $ when $n=1$.

For such a simplicial space, we can glue the disjoint union of the $\G_n \times \Delta ^n $ as follows, where $\Delta ^n$ is
the topological $n$-simplex. Let $ \delta_i : \Delta ^{n-1} \rightarrow \Delta ^n $ be the linear embedding of $\Delta ^{n-1}$
into $\Delta ^n$ as the $i$-th face. We define the classifying space of $\G$ as the geometric realization of its nerve as a 
simplicial space, that is:
\[ B \G = \mathop {\coprod } _n (\G_n \times \Delta ^n)/(d_i(g),x)\sim (g,\delta _i(x)). \]

This is usually called the fat realization of the nerve, meaning that we have chosen to leave out identifications involving
degeneracies. The two definitions will produce homotopy equivalent spaces provided that the topological category has sufficiently
nice properties. Another nice property of the fat realization is that if every $\G_n$ has the homotopy type of a CW-complex, then
the realization will also have the homotopy type of a CW-complex.

A strict homomorphism of groupoids $ \phi : \H \rightarrow \G $ induces a continuous map $B \phi : B \H \rightarrow B \G $. In 
particular, an important basic property is that an equivalence of groupoids induces a homotopy equivalence between classifying
spaces. This follows from the fact that an equivalence induces an equivalence of categories.
\newline

For any group $G$ we can construct the universal $G$-space $EG$ in the sense of Milnor or Milgram. We have analogous constructions
for a groupoid. The first construction is the analogue to the universal space of Milgram and it is also described in
Example 2.36 of \cite{GH}. Given the groupoid $\G $, construct the translation groupoid $\bar{\G} = \G \rtimes \G_1 $. This is the groupoid 
that has $\G_1$ as its space 
of objects and $ \G_1 \times _s \G_1 $ as its space of arrows, that is, only one arrow between two elements if they have the same source and 
none otherwise. The nerve of this category is given by:
\[ N\bar{\G }_k = \{ (g_1,\ldots,g_{k+1}) \mid g_i \in \G_1 \text{ for } i=1, \ldots, n \text{, and } s(g_1) = \cdots =
s(g_{k+1}) \}. \]
We consider the natural simplicial structure on $N\bar{\G}_*$ as the nerve of the category $\bar{\G}$. There is a natural action of $\G $ on 
$N\bar{\G }_k$ with anchor map $\pi $ given by the source map and $ (f_1,...,f_{k+1}) \cdot h = (f_1h,...,f_{k+1}h)$. With this action there are homeomorphisms $ N\bar{\G }_n / \G \cong \G_n $ for all $n>0$ given by the maps: 
\newline
\newline
\noindent $ \Phi _0: N\bar{\G}_0 / \G \to \G_0 $
\newline
\noindent $ [g] \mapsto t(g) $
\newline
\newline
\noindent $ \Phi _{2k+1}: N\bar{\G}_{2k+1} / \G \to \G_{2k+1} $ 
\newline
\noindent$ [g_1,\ldots,g_{2k+1}] \mapsto (g_2g_1^{-1},g_1g_3^{-1},g_4g_1^{-1}, \ldots, g_{2k+1}g_1^{-1}) $ 
\newline
\newline
\noindent$ \Phi _{2k}: N\bar{\G}_{2k} / \G \to \G_{2k} $ (if $k>0$)
\newline
\noindent$ [g_1,\ldots,g_{2k}] \mapsto (g_2g_1^{-1},g_1g_3^{-1},g_4g_1^{-1}, \ldots, g_1g_{2k}^{-1}) $ 
\newline
\newline
\noindent$ \Psi _{2k+1}: \G_{2k+1} \to N\bar{\G}_{2k+1} / \G $ 
\newline
\noindent$ (h_1,\ldots,h_{2k+1}) \mapsto [u(s(h_1)),h_1,h_2^{-1},h_3,\ldots,h_{2k+1}] $
\newline
\newline
\noindent$ \Psi _{2k}: \G_{2k} \to N\bar{\G}_{2k} / \G $ 
\newline
\noindent$ (h_1,\ldots,h_{2k}) \mapsto [u(s(h_1)),h_1,h_2^{-1},h_3,\ldots,h_{2k}^{-1}] $,
\newline
\newline
\noindent where $[g_1,\ldots,g_n]$ denotes the orbit of the element $(g_1,\ldots,g_n)$ in $ N\bar{\G }_n / \G $. These maps commute with the face 
operators $d_i$ for the simplicial structures on $N\bar{\G}_*$ and $\G_*$, therefore $ N\bar{\G } / \G $ and $ N\G $ are homeomorphic and so 
$ B\bar{\G } / \G \cong B\G $. It is also clear that $\G $ acts freely on $ B\bar{\G }$. 
\newline

The second construction imitates Milnor's universal $G$-space. Consider the spaces:
\[ E^n\G = \Big{\{} \mathop{\Sigma }_{i=1}^n \lambda _i g_i \mid s(g_1) = \cdots = s(g_n)\text{, } 
\mathop{\Sigma }_{i=1}^n \lambda _i = 1  \Big{\}} \]
\noindent with the subspace topology induced from the natural inclusion in the join of $n$ copies of $\G _1$ and an action of $\G$ given by the maps 
$ \displaystyle \pi \Big{(} \mathop{\Sigma }_{i=1}^n \lambda _i g_i \Big{)} = s(g_1) $ and
$ \displaystyle \Big{(} \mathop{\Sigma }_{i=1}^n \lambda _i g_i \Big{)} \cdot g = \mathop{\Sigma }_{i=1}^n \lambda _i g_ig $.

There are inclusions of $ E^n\G $ in $ E^{n+1}\G$ given by sending $ \displaystyle \Sigma _{i=1}^n \lambda _i g_i $ to 
$ \displaystyle \Sigma _{i=1}^{n+1} \mu _i h_i $, where $ h_i = g_i $, $ \lambda _i = \mu _i $ if $ i \leq n $ and $ h_{n+1} = u(s(g_1)) $, $ \mu_{n+1} = 0 $. Now 
define $ \displaystyle E\G = \lim _{\to} E^n\G $. Note that $\G $ acts freely on $ E^n\G $ for all $n$ and thus on $ E\G $. Let us denote $ B^n\G = E^n\G / \G $. Since 
the inclusions of of $ E^n\G $ in $ E^{n+1}\G$ described above are equivariant, they induce inclusions of $ B^n\G $ in $ B^{n+1}\G$.

The fibrant replacement of a topological groupoid (\cite{GH}, Definition 2.45) is weakly equivalent to
the original groupoid. It turns out that $ E\G $ is the geometric realization of the fibrant replacement of $\bar{\G}$ (\cite{GH}, Example 2.47) and 
therefore $E\G/\G $ and $ B\bar{\G} /\G = B\G$ are homotopy equivalent. From this point on we will use $E\G$ as our universal $\G$-space and
identify $E\G/\G $ with $B\G$.

\begin{defn}

Given a $\G $-space $X$, we define the Borel construction $ X_{\G } = ( X \times _{\pi } E\G ) / \G $.

\end{defn}

\begin{example}

Let $\G = H \rtimes M $. Then a $\G $-space is a $H$-space $X$ with a $H$-equivariant map to $M$. In this case, $\G $-equivariant 
vector bundles on $X$ correspond to $H$-vector bundles on $X$ and so $ K_{\G } (X) = K_H (X) $. We have $ E\G = M \times EH $ and 
$ B\G = M_H $.

\end{example}

\begin{example}

Let $M$ be a $\G $-space. Consider the groupoid $ \H = \G \rtimes M $. An $\H $-space is a $\G $-space $X$ with a 
$\G $-equivariant map to $M$. As in the previous example, $\H $-equivariant vector bundles on $X$ are just $\G $-equivariant 
vector bundles on $X$ and so $ K_{\H }(X) = K_{\G }(X) $. It can be shown that $ E\H = M \times _{\pi } E\G $ and $ B\H = M_{\G } $.

\end{example}

\begin{lemma} \label{FreeAction}

Let $X$ be a $\G$-space and $ Y = X \times _{\pi} E^n\G $. If $ Y/\G $ is compact, then 
there is an isomorphism $ K_{\G } (Y) \cong K ( Y / \G ) $.

\end{lemma}

\begin{proof}

Pulling back via the projection $ p_1: Y \to Y / \G $ takes vector bundles over $ Y / \G $ to $\G$-vector bundles over $Y$. On the other hand, if $E$ 
is a $\G$-vector bundle over $Y$, we claim that $ E / \G $ is a vector bundle over $Y / \G$. To prove this claim, we only need to show that 
$ E / \G $ is locally trivial. And to show this, it suffices to prove that $p_1$ has local sections around any point. 

Indeed, let $ [x,\mathop{\Sigma }_{i=1}^n \lambda _i g_i] = p_1(x,\mathop{\Sigma }_{i=1}^n \lambda _i g_i)$, where $x \in X $, $g_i \in \G _1$ and
$ \pi_X(x) = s(g_i) $ for $ i = 1, \ldots, n $. There is some $j$ for which $\lambda _j $ is not zero. Let $U$ be a $\G$-invariant open 
neighbourhood of $\mathop{\Sigma }_{i=1}^n \lambda _i g_i$ in $E^n\G$ such that every element of $U$ satisfies that $ \lambda _j $ is not zero, and
consider $ V = p_1(X \times _{\pi} U) $. Now consider the map $ f : V \to X \times_{\pi} U $ that sends $ [y,\mathop{\Sigma }_{i=1}^n \mu _i h_i] $ to 
$ (yg_j^{-1}, \mathop{\Sigma }_{i=1}^n \mu _i k_i) $, where $ k_i = h_i g_j^{-1} $ if $ i \neq j $ and $ k_j = ut(g_j)$. This map is well defined
and it is a section of $p_1$.

These two maps define a bijection between vector bundles over $Y/\G$ and $\G$-vector bundles on $Y$. Hence any $\G$-vector bundle on $Y$ is the pullback
of a vector bundle on $Y / \G $. To prove the lemma we need to show that vector bundles on $ Y/\G$ pull back to extendable $\G$-vector bundles on $Y$. 
The anchor map $ \pi _1 : Y \longrightarrow \G_0 $ induces a map $ \pi _2 : Y/\G \longrightarrow | \G | $. These maps fit into a commutative 
diagram:
\[
\diagram
Y \rto^{\pi_1} \dto^{p_1} & \G_0 \dto^{p_2} \cr
Y/\G \rto^{\pi _2} & |\G|,
\enddiagram
\]
where $p_2$ is the projection given in Definition \ref{Orbits}. 

Let $r$ and $q$ be the only maps from $Y/\G$ and $|\G|$ to a point, respectively. Given a vector bundle $A$ over $Y/\G$, there is a vector bundle 
$V$ over a point such that $r^*V = A \oplus B $ for some vector bundle $B$ over $Y/\G$. If $W=q^*V$, then we also have $ \pi _2^*W = A \oplus B $. 
Consider $W' = p_2^*W$. We have
\[ \pi_1 ^*W' = \pi_1 ^*p_2^*W = p_1^*\pi _2^*W = p_1^* (A \oplus B) = p_1^*A \oplus p_1^*B.   \qedhere \]
\end{proof}

Note that if $ B^n\G$ is compact, we obtain $ K_{\G } (E^n\G) \cong K ( B^n\G) $ by taking $ X = \G_0 $. In the next lemma and in the rest of the paper 
we will use representable $K$-theory for spaces that are not compact.

\begin{lemma} \label{Nilpotence}

Consider $ \displaystyle \tilde{M} = \Big{\{} \mathop{\Sigma }_{j=1}^n \lambda _j g_j \mid g_1 = \cdots = g_n\text{, } 
\lambda _1 = \cdots = \lambda _n = 1/n  \Big{\}} $, which is a $\G$-equivariant subspace of $ E^n\G $ and $ M = \tilde{M}/\G \subset B^n\G $. 
Then any product of $n$ elements in $ K^* (B^n\G , M) $ is zero. 

\end{lemma}

\begin{proof}

The result is obvious if $ n = 1 $. Let $ n > 1 $ and consider the $\G$-subspaces   
$ \tilde{A}_i = \Big{\{} \mathop{\Sigma }_{j=1}^n \lambda _j g_j \mid \lambda _i \geq 1/n \Big{\}} $ of $ E^n\G $ and the corresponding
subspaces $ A_i = \tilde{A}_i / \G $ of $ B^n\G $. 
Given $h \in \G _1$, let us denote by $h(i)$ the element $ \mathop{\Sigma }_{j=1}^n \lambda _j g_j \in E^n\G $ with $ \lambda_i = 1 $, $ \lambda_j = 0 $
if $ j \neq i $ and $ g_j = h $ for all $j$. Now, we define $ \displaystyle \tilde{M_i} = \{ h(i) \mid h \in \G _1 \} $. The space $\tilde{M_i}$ is a deformation 
retract of the spaces $ \tilde{A_i} $, as the following maps show:
\newline
\newline
\noindent $ \rho _i: \tilde{A}_i \to \tilde{M}_i $
\newline
\noindent $ \displaystyle \rho_i \left( \mathop{\Sigma }_{j=1}^n \lambda_j g_j \right) = g_i(i) $ 
\newline
\newline
\noindent $ H_i: \tilde{A}_i \times I \to \tilde{A}_i $
\newline
\noindent $ \displaystyle H_i\left( \mathop{\Sigma }_{j=1}^n \lambda_j g_j,t \right) = \mathop{\Sigma }_{j=1}^n \mu_j h_j $, where $ \mu_i = t\lambda_i + (1-t) $,
$\mu _j = t\lambda_j$ if $ j \neq i $ and $ h_j = g_i $ for all $j$.
\newline
\newline
Since all these maps are $\G$-equivariant, $M_i = \tilde{M}_i / \G $ is a deformation retract of the space $A_i$. 
In particular $ K^*(A_i,M_i) = 0 $. The map $ \tilde{\Phi} : \tilde{M}_i \to \tilde{M} $ given by 
$ \Phi ( g(i)) = \mathop{\Sigma }_{j=1}^n \frac{1}{n} g $ is a homeomorphism and it is $\G$-equivariant, hence it
defines a homeomorphism $ \Phi : M_i \to M $. Then the relative $K$-theory of the pair $(A_i,M)$ is isomorphic to 
the $K$-theory of the mapping cone for the map $ \Psi : M_i \to A_i $ that takes $ g(i) $ to 
$ \mathop{\Sigma }_{j=1}^n \frac{1}{n} g \in A_i $. But the inclusion of $ M_i $ in $ A_i $ is homotopic to $ \Psi $ via the
homotopy   
\newline
\newline
\noindent $ H_i: \tilde{M}_i \times I \to \tilde{A}_i $
\newline
\noindent $ \displaystyle H_i\left( g(i),t \right) = \mathop{\Sigma }_{j=1}^n \lambda_j g $, where $ \displaystyle \lambda_i = t + \frac{1-t}{n} $ 
and $ \displaystyle \lambda_j = \frac{1-t}{n} $ if $ j \neq i $.
\newline
\newline
Therefore $ K^*(A_i,M) \cong K^*(A_i,M_i) = 0 $, and so $ K^*(B^n\G,M) \cong K^*(B^n\G,A_i) $ for all $ i $ via the long exact sequence for the triples 
$(B^n\G,A_i,M)$. Since $ B^n\G = A_1 \cup A_2 \cup \ldots \cup A_n $, the relative product 
$ K^*(B^n\G,A_1) \otimes K^*(B^n\G,A_2) \otimes \cdots \otimes K^*(B^n\G,A_n) \to K^*(B^n\G,A_1 \cup A_2 \cup \ldots \cup A_n) $ is
identically zero, and therefore the product of $n$ elements of $ K^*(B^n\G,M) $ is zero.
%
\end{proof}

\begin{defn}

A Lie groupoid $\G$ is finite if $\G_0$ is a finite $\G$-CW-complex and the spaces $B^n\G$ are compact.

\end{defn}

Let $\G$ be a finite groupoid. We have an augmentation homomorphism $ K_{\G}^*(\G_0) \rightarrow K^*(\G_0) $ given by forgetting the $\G $-action. 
Let $ I_{\G } $ be the kernel of this map and let us consider the following composition:
\[ K_{\G }^* (\G_0) \stackrel{\pi ^*}{\rightarrow} K_{\G }^* (E^n\G ) \stackrel{\cong}{\rightarrow} K^*(B^n\G ) \stackrel{i^*}{\rightarrow} 
K^*(M) \stackrel{\cong}{\rightarrow} K^*(\G_0), \]
where $ \pi : E^n\G \to \G_0 $ is the anchor map, the middle isomorphism is given by Lemma \ref{FreeAction}, 
$ i: M \to B^n\G $ is the inclusion given in Lemma \ref{Nilpotence} and the last isomorphism is given by the homeomorphism $ M \to \G _0$ 
induced by the $\G$-homeomorphism $ \tilde{M} \to \G _1 $ that sends $ \mathop{\Sigma }_{j=1}^n \frac{1}{n} g $ to $ g $. The composition of 
these four maps coincides with the augmentation homomorphism, whose whose kernel is $ I_{\G } $, so the map 
$K_{\G }^* (\G_0) \rightarrow K_{\G }^* (E^n\G )$ factors through $ K_{\G }^* (\G_0)  / I_{\G }^n $.

For any $\G $-space $X$, $K_{\G }^* (X) $ is a module over $ K_{\G }^* (\G_0) $ and by naturality the homomorphism 
$ K_{\G }^* (X) \rightarrow K_{\G }^* ( X \times _{\pi } E^n\G ) $ factors through $ K_{\G }^* (X) / I_{\G }^n K_{\G }^* (X) $, giving a 
homomorphism $ \Phi _n: K_{\G }^* (X) / I_{\G }^n K_{\G }^* (X) \rightarrow K_{\G }^* (X \times _{\pi } E^n\G ) $. The equivariant inclusion
$ E^n\G \to E^{n+1}\G $ and $ I_{\G}^{n+1} \subset I_{\G}^n $ induce the following commutative diagram:
\[
\diagram
K_{\G }^* (X) / I_{\G }^{n+1} K_{\G }^* (X) \dto \rto^{\Phi_{n+1}} & K_{\G }^* (X \times _{\pi } E^{n+1}\G ) \dto \cr
K_{\G }^* (X) / I_{\G }^n K_{\G }^* (X) \rto^{\Phi_n} & K_{\G }^* (X \times _{\pi } E^n\G ),
\enddiagram
\]
which induces a homomorphism of pro-rings $ \{ K_{\G }^* (X) / I_{\G }^n K_{\G }^* (X) \} \rightarrow \{ K_{\G }^* (X \times _{\pi } E^n\G ) \} $.
\begin{conj}

Let $\G$ be a finite Lie groupoid and $X$ a $\G$-space. Then we have an isomorphism of pro-rings 
\[ \{ K_{\G}^*(X)/I_{\G}^n K_{\G}^*(X) \} \longrightarrow \{ K_{\G}^*(X \times _{\pi} E^n\G ) \}. \]
\end{conj}

If a groupoid $\G$ satisfies this conjecture for $ X = \G_0 $, we will say $\G$ satisfies the completion theorem.

\section{The completion theorem}

Recall that we are using representable $K$-theory for spaces that are not compact.

\begin{lemma} \label{ClassicalCompletion}

Let $\G = G \rtimes X $, where $G$ is a compact Lie group and $X$ is a compact $G$-space such that $K_G^*(X)$ is finitely
generated over $R(G)$. Then $\G$ satisfies the completion theorem. 

\end{lemma}

\begin{proof}

Let $I_X$ be the kernel of $ K_G^*(X) \longrightarrow K^*(X) $. We would like to prove that there is an isomorphism of 
pro-groups $ \{ K_G^*(X) / I_X^n \} \cong \{ K_G^*(X \times E^nG) \} $. By the Atiyah-Segal completion theorem (Theorem 2.1 in \cite{AS}), 
we have an isomorphism $ \{ K_G^*(X) / I_G^nK_G^*(X) \} \cong \{ K_G^*(X \times E^nG) \} $. So it suffices to prove that the $I_G$-adic topology
and the $I_X$-adic topology are the same in $K_G^*(X)$. Since $K_G^*(X)$ is a module over $R(G)$, we have
$ I_GK_G^*(X) \subset I_X $.

Let $K_n$ be the kernel of $\alpha _n : K_G^*(X) \longrightarrow K_G^*(X \times E^nG) $. By Corollary 2.3 in \cite{AS}, 
the sequence of ideals $ \{ K_n \} $ defines the $I_G$-adic topology on $K_G^*(X)$. In
particular, there is $m \in \N $ such that $ K_m \subset I_GK_G^*(X) $. Note that $ K_1 = I_X $. Consider the composition
$ K_G^*(X) \longrightarrow K_G^*(X \times E^mG) \longrightarrow K_G^*(X \times E^1G) \cong K^*(X) $. Since $ X \times E^mG $ is
the union of $m$ open sets which are homotopy equivalent to $ X \times E^1G = X \times G $, we have 
$K_G^*(X \times E^mG,X \times E^1G)^m = 0 $. Thus the first map factors through $I_X^m$, thus $ I_X^m \subset K_m $. Hence
$I_X^m \subset I_GK_G^*(X) $.
\end{proof}

\begin{lemma}

Let $ \Theta : \H \to \G$ be a local equivalence of finite Lie groupoids. Then $\G$ satisfies the completion theorem if and only if $\H$ does.

\end{lemma}

\begin{proof} The local equivalence $ \Theta : \H \longrightarrow \G $ induces an 
isomorphism $ f : K_{\G}^*(\G_0) \stackrel{\cong}{\longrightarrow} K_{\H}^*(\H_0) $. The following diagram is commutative: 
\[
\diagram
K_{\G}^*(\G_0) \rto^{f} \dto^{\alpha } & K_{\H}^*(\H_0) \dto^{\beta } \cr
K^*(\G_0) \rto^{\phi} & K^*(\H_0).
\enddiagram
\]
Therefore we have $  f(I_{\G}) \subset I_{\H} $. Let $g$ be the inverse of $f$, $x \in I_{\H}$ and $y = g(x) $. Since $\beta (x) =0$,
we have $\phi \alpha (y) = 0 $. But $\alpha(y) = ( n_y,a_y) \in \Z \oplus \tilde{K}^*(\G_0) \cong K^*(\G_0) $ and so 
$ \phi \alpha (y) = ( n_y , \tilde{ \phi } (a_y)) $. This implies $n_y=0$, that is, $ \alpha g(I_{\H}) \subset \tilde{K}^*(\G_0) $. 

Since $\G_0$ is compact, there is $ m \in \N $ such that $ \tilde{K}^*(\G_0)^m = 0 $. Then we have $\alpha g(I_{\H}^m) \subset \tilde{K}^*(\G_0)^m = 0 $ 
and so $ g(I_{\H}^m) \subset I_{\G} $. Thus $ I_{\H}^m = fg(I_{\H}^m) \subset f(I_{\G}) $ and the topologies induced by $I_{\G}$
and $I_{\H}$ on $ K_{\H}^*(\H_0) $ are the same. Therefore we have an isomorphism of pro-rings 
$ \{ K^*(\G_0)/I_{\G}^n \} \cong \{ K^*(\H_0)/I_{\H}^n \} $.

The local equivalence $\Theta $ also induces a homotopy equivalence between $B\G$ and $B\H$. If we consider the associated filtrations to each of
these spaces, $\{ B^n\G \} $ and $ \{ B^n\H \} $, the local equivalence must take $ B^n\G$ to some $ B^{n+k}\H $ and $B^n\H$ to some
$ B^{n+k'}\G $. Hence we have an isomorphism of pro-rings $ \{ K^*(B^n\G) \} \cong \{ K^*(B^n\H) \} $. The lemma follows
from the diagram:
\[
\diagram
\{ K^*(\G_0)/I_{\G}^n \} \rto^{\cong} \dto & \{ K^*(\H_0)/I_{\H}^n \} \dto \cr
\{ K^*(B^n\G) \} \rto^{\cong} & \{ K^*(B^n\H) \} \qedhere
\enddiagram  
\] 
\end{proof}

From the previous lemma, we obtain the following theorem:

\begin{thm} \label{WeakEquivalence}

Let $\G$ and $\H$ be weakly equivalent finite groupoids. Then $\G$ satisfies the completion theorem if and only if 
$\H$ does.

\end{thm}

Note that Theorem \ref{WeakEquivalence} and Lemma \ref{ClassicalCompletion} imply the following:

\begin{cor}

If $\G$ is Bredon-compatible and $U$ is a $\G$-cell, $K_{\G}^*(U) $ is a finitely generated abelian group and the groupoid
$\G \rtimes U $ satisfies the completion theorem.

\end{cor}

This corollary tells us that the completion theorem is true for $\G$-cells. Now we will use this corollary
to prove the completion theorem for finite $\G$-CW-complexes. We will also need to use the following lemma:

\begin{lemma}

Fix any commutative Noetherian ring $A$, and any ideal $ I \subset A $. Then for any exact sequence $ M' \rightarrow M \rightarrow M''$
of finitely generated $A$-modules, the sequence
\[ \{ M'/I^nM' \} \longrightarrow \{ M / I^nM \} \longrightarrow \{ M''/I^nM'' \} \]
of pro-groups (pro-$A$-modules) is exact.

\end{lemma}

\begin{proof}

See \cite{LO}, Lemma 4.1.
\end{proof}

Let $X$ be a finite $\G$-CW-complex with $p$-cells $ \{ V_i = U_i \times D^p \mid i \in I_p \} $ (see Definition
\ref{B}) and consider the spectral sequence for this decomposition in $\G $-equivariant $K$-theory:
\[  E_1 ^{pq} = \mathop{\mathop{\prod}}   _{i \in I_p} K_{\G }^q (V_i) \Longrightarrow K_{\G }^{p+q} (X). \]
%
%
This is a spectral sequence of $ K_{\G}^*(\G_0) $-modules. Let us assume that $\G$ is finite so that
$ K_{\G}^*(\G_0) $ is a Noetherian ring. All elements in these spectral sequences are finitely generated, and
so by the previous lemma, the functor taking a $ K_{\G}^*(\G_0) $-module $M$ to the pro-group 
$\{ M/I_{\G}^nM \} $ is exact and we can form the following spectral sequence of pro-rings: 
\[  F_1 ^{pq} = \left\{ \mathop{\mathop{\prod}}   _{i \in I_p} K_{\G }^q (V_i)/ I_{\G }^nK_{\G }^q
(V_i) \right\}
\Longrightarrow \{ K_{\G }^{p+q} (X)/ I_{\G }^nK_{\G }^{p+q} (X) \}. \]

Now consider the $\G$-maps $ h_n : X \times _{\pi} E^n\G  \longrightarrow X $. They give us another spectral sequence of
pro-rings:
\[  \bar{F}_1 ^{pq} = \left\{ \mathop{\mathop{\prod}}   _{i \in I_p} K_{\G }^q (h_n^{-1}V_i) \right\}
\Longrightarrow \{ K_{\G }^{p+q} (X \times _{\pi} E^n\G) \}. \]
We also have a map of spectral sequences $ \phi : F \longrightarrow \bar{F} $ induced by the projections
$ h_n^{-1}V_i = V_i \times _{\pi} E^n/G \to V_i $.

If $\G$ is Bredon-compatible, the groupoids $ \G \rtimes V_i $ satisfy the completion theorem for all $i$. Since we are 
taking quotient by the ideal $I_{\G}$ and not by $I_{\G \rtimes V_i}$, we need to check that both topologies
coincide. We consider the long exact sequences in equivariant and non-equivariant $K$-theory for the pair $(C_V,V)$ where $V$ is any $V_i$ and $C_V$ is the mapping cylinder
of the map $ \pi : V \longrightarrow \G_0 $. Note that $C_V$ is $\G$-homotopy equivalent to $\G_0$ and $ V
\subset C_V $.
\[
\diagram
K_{\G}(C_V,V) \rto \dto & K_{\G}(\G_0) \rto \dto & K_{\G}(V) \rto \dto & K_{\G}^1(C_V,V) \dto \cr
K(C_V,V) \rto & K(\G_0) \rto & K(V) \rto & K^1(C_V,V).
\enddiagram
\]
It is clear that $I_{\G}K_{\G}(V) \subset I_V $. Now let $m \in \N$ such that $K(C_V,V)^m=0$ and $ n \in \N $ such that
$K_{\G}(C_V,V)^n = 0 $. Then $I_V^{nm} \subset I_{\G}K_{\G}(V) $ and so the topologies coincide.

This proves $\phi$ is an isomorphism when restricted to any particular element $F^{ij}$ and therefore, it is an isomorphism
of spectral sequences. In particular, we have an isomorphism $ \{ K_{\G }^{p+q} (X)/ I_{\G }^nK_{\G }^{p+q} (X) \} \cong 
\{ K_{\G }^{p+q} (X \times _{\pi} E^n\G) \} $. To summarize, we have proved:

\begin{thm} \label{Completion}
Let $\G$ be a Bredon-compatible finite Lie groupoid and $X$ a finite $\G$-CW-complex. Then we have an isomorphism of pro-rings 
\[ \{ K_{\G}^*(X)/I_{\G}^n K_{\G}^*(X) \} \longrightarrow \{ K^*_{\G}(X \times _{\pi} E^n\G ) \}. \]
\end{thm}

\begin{cor}

Under the same circumstances, the homomorphism $ K_{\G}^*(X) \rightarrow K^*(X_{\G}) $ induces an isomorphism of the $I_{\G}$-adic 
completion of $ K_{\G}^*(X) $ with $K^*(X_{\G}) $.

\end{cor}

\section{Applications}

Unless otherwise stated, throughout this whole section $S$ will be a Lie group, but not necessarily compact. To study proper actions of $S$, we
can consider the groupoid $\G = S \rtimes \underline{E}S$, where $\underline{E}S$ is the universal space for proper actions of 
$S$ (\cite{L}, Definition 1.8). This space is a proper $S$-CW-complex such that $\underline{E}S^{G} $ is contractible for all 
compact Lie subgroups $G$ of $S$ and such that every proper $S$-CW-complex has an $S$-map to $\underline{E}S$ that is unique up to $S$-homotopy. 
The existence of $\underline{E}S$ is shown in Theorem 1.9 of \cite{L}. Some immediate consequences follow:

\begin{itemize}

\item Proper $S$-CW-complexes are $\G$-CW-complexes.

\item Extendable $\G$-bundles on a proper $S$-CW-complex $X$ are extendable $S$-bundles for any $S$-map 
$ X \longrightarrow \underline{E}S $, since all of them are $S$-homotopic.

\end{itemize}

In order for $\G $ to be finite, we need $\underline{E}S $ to be a compact finite $S$-CW-complex and the spaces $B^n\G$ to
be compact. Note that this is not necessarily the case, but when these conditions hold we are under the hypotheses of Theorem \ref{Completion}.

Equivariant $K$-theory for actions of finite groups and compact Lie groups have been studied extensively in the past. It is
a well-known fact that for these actions, $\G $ is Bredon-compatible \cite{SG3} and finite and so they satisfy a completion theorem, first 
proved in \cite{AS} using different methods. 

Proper actions of discrete groups and totally disconnected groups that are projective limits of discrete groups also satisfy Bredon-compatibility, 
as shown in \cite{LO} and \cite{S}, and therefore a completion theorem follows if the finiteness condition holds. A completion theorem for
proper actions of discrete groups is also proved in \cite{LO}.

In general, vector bundles may not be enough to construct an interesting equivariant cohomology theory for proper actions of
second countable locally compact groups \cite{Ph}, but they suffice for two important families, almost compact groups and matrix
groups \cite{Ph2}.

Almost compact groups, that is, second countable locally compact groups whose group of connected components are 
compact, always have a maximal compact subgroup. Any space with a proper action of one of these groups is the induction of a
space with an action of that compact subgroup and so the study of proper actions of almost compact groups are reduced to studying
compact Lie group actions. This is carried out in \cite{Ph2}, where it is proven that these action groupoids are Bredon-compatible
that a completion theorem holds. This is done by showing that the completion maps are compatible with the
reduction map to the maximal compact subgroup. With different techniques it is proven that proper actions of matrix groups, that is,
closed subgroups of $GL(n,\R)$, are Bredon-compatible. Proper actions of abelian Lie groups constitute a particular instance of 
this case.

The groupoid $\G$ is not necessarily Bredon-compatible. When $ \G $ is a Bredon-compatible groupoid we must have
Vect$_{\G}(S/G) = $ Vect$_G(pt) $ for a compact subgroup $G$ of $S$. Let $S$ be a Kac-Moody group and $T$ its maximal 
torus. Note that $S$ is not a Lie group, but all the constructions and results in Sections 3 and 4 remain valid for a 
topological groupoid with averaging of vector bundles.

There is an $S$-map $ S/T \longrightarrow \underline{E}S $ which is unique, up to homotopy. Given an $S$-vector bundle $V$ on 
$\underline{E}S$, the pullback to $S/T$ is given by a finite-dimensional representation of $T$ invariant under the Weyl group.
This representation gives rise to a finite-dimensional representation of $S$. All finite-dimensional representations of
Kac-Moody groups which are not compact Lie groups are trivial. Therefore $V$ must be trivial and so extendable $S$-vector 
bundles on $S/T$ only come from trivial representations of $T$. 

In order to deal with these groups, it is more convenient to use dominant $K$-theory, which was developed in \cite{K}. Kac-Moody
groups possess an important class of representations called dominant representations. A dominant representation of a Kac-Moody 
group in a Hilbert space is one that decomposes into a sum of highest weight representations. Equivariant 
$K$-theory for proper actions of Kac-Moody groups is defined as the representable equivariant cohomology theory modeled on
the space of Fredholm operators on a Hilbert space which is a maximal dominant representation of the group.

\end{document}